\documentclass[spmpsci]{sn-jnl}

\usepackage{subcaption}
\usepackage{xcolor}
\usepackage{enumerate}

\jyear{2021}%

\theoremstyle{thmstyleone}%
\newtheorem{theorem}{Theorem}

\theoremstyle{thmstyletwo}%
\newtheorem{example}{Example}%
\newtheorem{remark}{Remark}%

\theoremstyle{thmstylethree}%
\newtheorem{definition}{Definition}%
\newtheorem{proposition}{Proposition}%
\newtheorem{lemma}{Lemma}%
\newtheorem{corollary}{Corollary}%
\begin{document}

\title[On the metric representation of the vertices of a graph]{On the metric representation of the vertices of a graph}


\author[1]{\fnm{Merc\`e} \sur{Mora}}\email{merce.mora@upc.edu}

\author*[2,3]{\fnm{Mar\'ia Luz} \sur{Puertas}}\email{mpuertas@ual.es}

\affil[1]{\orgdiv{Department of Mathematics}, \orgname{Universitat Polit\`ecnica de Catalunya}, \orgaddress{\street{Carrer Jordi Girona, 31}, \city{Barcelona}, \postcode{08034}, \country{Spain}}}

\affil[2]{\orgdiv{Department of Mathematics}, \orgname{Universidad de Almer\'ia}, \orgaddress{\street{Carretera Sacramento s/n}, \city{Almer\'ia}, \postcode{04120}, \country{Spain}}}

\affil[3]{\orgdiv{Agrifood Campus of International Excellence (ceiA3)}, \orgname{Universidad de Almería}, \orgaddress{\street{Carretera Sacramento s/n}, \city{Almer\'ia}, \postcode{04120}, \country{Spain}}}


\abstract{The metric representation of a vertex $u$ in a connected graph $G$ respect to an ordered vertex subset $W=\{\omega_1, \dots , \omega_n\}\subset V(G)$ is the vector of distances $r(u\vert W)=(d(u,\omega_1), \dots , d(u,\omega_n))$. A vertex subset $W$ is a resolving set of $G$ if $r(u\vert W)\neq r(v\vert W)$, for every  $u,v\in V(G)$ with $u\neq v$. Thus, a resolving set with $n$ elements provides a set of metric representation vectors $S\subset \mathbb{Z}^n$ with cardinal equal to the order of the graph. In this paper, we address the reverse point of view, that is, we characterize the finite subsets $S\subset \mathbb{Z}^n$ that are realizable as the set of metric representation vectors of a graph $G$ with respect to some resolving set $W$. We also explore the role that the strong product of paths plays in this context. Moreover, in the case $n=2$, we characterize the sets $S\subset \mathbb{Z}^2$ that are uniquely realizable as the set of metric representation vectors of a graph $G$ with respect to a resolving set $W$.
}




\maketitle

\section{Introduction}\label{sec:intro}

The metric dimension of a connected graph is a well known parameter introduced in the seventies~\cite{Harary1976,Slater1975}. Since then, it has been extensively studied from many different points of view, such as variations of the original definition, computational aspects and practical applications~\cite{Khuller1996}. Two recent surveys about this topic~\cite{Yero2022,Tillquist2021} provide an extensive review of the research developed around it.

All the graphs considered in  this paper are finite, simple  and connected. We quote the definition of the metric dimension from~\cite{Chartrand2000}, which emphasizes the role of the coordinates.

\begin{definition}[\cite{Chartrand2000,Harary1976,Slater1975}]
The \emph{metric representation} of a vertex $u\in V(G)$ with respect to a vertex subset $W=\{\omega_1, \dots , \omega_n\}\subset V(G)$ is $r(u\vert W)=(d(u,\omega_1), \dots d(u,\omega_n))$. A vertex subset $W$ of a graph $G$ is a \emph{resolving set} if $r(u\vert W)\neq r(v\vert W)$, for each pair of vertices $u,v\in V(G)$, with $u\neq v$. The \emph{metric dimension} of $G$, denoted by $\dim (G)$, is the minimum cardinality of a resolving set of $G$.
A \emph{metric basis} of $G$ is a resolving set of cardinality $\dim (G)$.\end{definition}

Hence, given a graph $G$ and a resolving set $W$ with $n$ vertices, a finite set of metric representation vectors $S=\{ r(u\vert W)\colon u\in V(G)\}\subset \mathbb{Z}^n$ is provided. Some results related to properties shared by the graphs representing the same vector set can be found in~\cite{Feit2016,Feit2018}. In this paper, we address the problem from a different point of view and we focus on the following questions. Given a finite subset $S\subset \mathbb{Z}^n$, is there a graph $G$ such that $S$ is the set of metric representation vectors for some resolving set $W$? If the answer is positive, is such graph unique? Moreover, we want to give the answers to both questions just by means of elementary properties of the elements of $S$.

The following examples show that, if $S\subset \mathbb{Z}^n$ is a set of metric representation vectors for some graph $G$ and some resolving set $W$, neither the resolving set $W$, nor the graph $G$ are necessarily unique. Moreover, there exist resolving sets of non-isomorphic graphs that give rise to the same metric representation set.

\begin{example}
The subset $S=\{(0,1), (0,2), (1,1)\}$ of $\mathbb{Z}^2$ is not a set of metric representation vectors of any graph. Suppose on the contrary that there exists a graph $G$ and a resolving set $W=\{\omega_1, \omega_2\}$ such that $S=\{r(u\vert W)\colon u\in V(G)\}$. Then, there exist $u,v\in V(G)$ satisfying $(0,1)=(d(u,\omega_1), d(u,\omega_2))$ and $(0,2)=(d(v,\omega _1),d(v,\omega_2))$. This means that $u=v=\omega _1$, but $d(\omega_1,\omega_2)=d(u,\omega_2)=1\neq 2=d(v,\omega _2)=d(\omega_1,\omega_2)$, a contradiction.
\end{example}

\begin{example}\label{ex:k5}
Consider the set $S=\{ (0,1,1,1), (1,0,1,1), (1,1,0,1),(1,1,1,0), (1,1,1,1)\}$  of $\mathbb{Z}^4$. Clearly, any subset $W$ with four vertices of a complete graph of order 5  satisfies  $S=\{r(u\vert W)\colon u\in V(K_5)\}$. Therefore, $S$ is the set of metric representation vectors of several resolving sets of the same graph.
\end{example}

\begin{example}
Let $S=\{ (0,2), (2,0), (1,1), (3,1), (1,3)\}\subset \mathbb{Z}^2$ and consider the path $P_5$ depicted in Figure~\ref{fig:path_5}. Then, $W=\{u_2,u_4\}$ is a resolving set of $P_5$ and $S=\{r(u\vert W)\colon u\in V(P_5)\}$. Moreover, it can be easily checked that $W$ is the only resolving set of $P_5$ with this metric representation vector set.

\begin{figure}[h]
    \centering
    \includegraphics[width=0.25\textwidth]{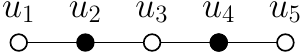}
    \caption{A resolving set of the path $P_5$.}
    \label{fig:path_5}
\end{figure}

\end{example}

\begin{example}\label{ex:non-isomorphic}
In Figure~\ref{fig:dimension_3}, two non isomorphic graphs $G_1$ and $G_2$ are shown. In both cases, the set $W=\{\omega_1, \omega_2, \omega_3\}$ (black vertices) is a resolving set. Moreover, $\{r(u\vert W)\colon u\in V(G_1)\} =\{r(v\vert W)\colon v\in V(G_2)\}=\{(0,1,2), (1,0,1), (2,1,0), (1,2,3), (2,3,2), (3,2,1), (2,1,2), (1,2,1)\}$. Therefore, two non isomorphic graphs can share the same metric representation vector set for some resolving sets.

\begin{figure}[ht]
     \centering
     \begin{subfigure}{.4\textwidth}
     \centering
         \includegraphics[width=\textwidth]{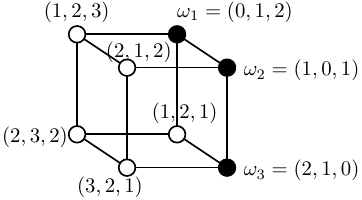}
         \caption{ $G_1$}
         \label{fig:dimension_3_a}
     \end{subfigure}
     \hspace{1cm}
     \begin{subfigure}{.4\textwidth}
     \centering
         \includegraphics[width=\textwidth]{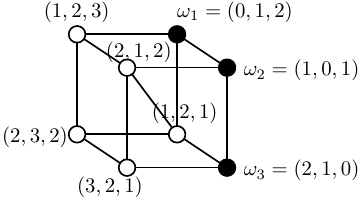}
         \caption{$G_2$}
         \label{fig:dimension_3_b}
     \end{subfigure}
        \caption{Two non isomorphic graphs sharing the same metric representation vector set.}
        \label{fig:dimension_3}
\end{figure}
\end{example}

The paper is organized as follows. In Section~\ref{Sec:2}, we introduce the concept of set realizable as the metric representation vectors of a graph and we characterize such sets by means of their internal features. In Section~\ref{Sec:3}, we study the problem of the uniqueness of the realization, introducing the notion of equivalent realizations. To this end, we first explore in depth the relationship between realizable sets and the strong product of paths. Finally, we present necessary and sufficient conditions to ensure the uniqueness of realizations for subsets of $\mathbb{Z}^2$.

From now on, $[n]$ denotes the set $\{1,2,\dots, n\}$, for every integer $n\ge 1$. Also, if $x\in \mathbb{Z}^n$, then $x_i$ denotes its $i$-th coordinate, that is, $x=(x_1,\dots ,x_n)$.
\section{Realizable sets}\label{sec:characterization}\label{Sec:2}

In this section we present some basic properties of the vector sets that are realizable as the metric representation vector set of a graph. Moreover, we characterize such vector sets in terms of some of those properties. First of all, we pose the following definition.

\begin{definition}
A subset $S\subset \mathbb{Z}^n$ is \emph{realizable} if there exists a graph $G$  and a resolving set $W$ of $G$ such that $S=\{r(u\vert W)\colon u\in V(G)\}$.
In such a case, we say that $(G,W)$ is a \emph{realization} of $S$.
\end{definition}

From now on, if $(G,W)$ is a realization of $S$, then for each $x\in S$, we denote by $u_G(x)$ the only vertex in $G$ such that its metric representation with respect to $W$ is $x$, so that $r(u_G(x)\vert W)=x$. We write simply $u(x)$  if the graph $G$ is clear from context.

The following properties can be easily derived from the definition.

\begin{proposition}\label{pro:basic}
Let $(G,W)$ be a realization of a set $S\subset \mathbb{Z}^n$. Then:
\begin{enumerate}
    \item $\vert V(G)\vert =\vert S\vert$.
    \item $\vert W\vert =n$ and therefore, $\dim (G)\leq n\leq \vert S\vert $.
    \item If $x\in S$, then $x_i\geq 0$, for every $i\in [n]$. Moreover, $x$ has at most one coordinate equal to zero.
    \item For every $i\in[n]$, there exists a unique $x\in S$ such that $x_i=0$.
\end{enumerate}
\end{proposition}

\begin{proof}
\begin{enumerate}
    \item We have $\vert V(G)\vert =\vert \{r(u\vert W)\colon u\in V(G)\}\vert =\vert S\vert$, since $W$ is a resolving set of $G$.
    \item If $v\in V(G)$, then $r(v\vert W)\in S\subset \mathbb{Z}^n$. Hence, $r(u\vert W)$ is a vector with $n$ coordinates and $\vert W\vert =n$. Therefore, $\dim(G)\leq \vert W\vert=n\leq \vert V(G)\vert =\vert S\vert$.
    \item Suppose that $W=\{\omega _1,\dots, \omega_n\}$. Then,  $x_i=d(u(x),\omega _i)\geq 0$,  for every $x\in S$ and for every $i\in [n]$. Moreover, if $x\in S$ satisfies $x_i=x_j=0$, then
    $u(x)=\omega_i=\omega_j$, that implies $i=j$. Therefore, $x$ has at most one coordinate equal to zero.
    \item Let $i\in[n]$ and take $x=r(\omega_i \vert W)$. Clearly $x\in S$ and  $x_i=d(\omega_i,\omega_i)=0$.
    Moreover, if $y\in S$ satisfies $y_i=0$, then $u(y)=\omega_i$, because
    $y_i=d(u(y),\omega_i)=0$. Therefore,  $x=y=r(\omega_i\vert W)$, as desired.
\end{enumerate}
\end{proof}

\begin{remark}\label{rem:dimension}
For every graph $G$ and every resolving set $W$ of $G$, it is obvious from the definition that $(G,W)$ is a realization of the set of distance vectors $\{ r(u\vert W)\colon u\in V(G)\}\subset \mathbb{Z}^{\vert W\vert}$.
\end{remark}

\begin{corollary}\label{cor:md_characterization}
Let $G$ be a graph and $n$ an integer such that $n\le \vert V(G) \vert$. Then, $\dim(G)\leq n$ if and only if there is a resolving set $W$ of $G$ such that $(G,W)$ is a realization of a set $S\subset \mathbb{Z}^n$.
\end{corollary}
\begin{proof}
If $(G,W)$ is a realization of $S\subset \mathbb{Z}^n$, by Proposition~\ref{pro:basic} we know that $\dim (G)\leq n$. Conversely, if $\dim(G)\leq n$ then, there exists a metric basis $M$ such that $\vert M\vert=\dim (G)$. Let $W\subset V(G)$ be such that $M\subset W$ and $\vert W\vert =n$. Clearly $W$ is a resolving set of $G$ and, by Remark~\ref{rem:dimension}, $(G,W)$ is a realization of $S=\{r(u\vert W)\colon u\in V(G)\}\subset \mathbb{Z}^{\vert W\vert}=\mathbb{Z}^n$.
\end{proof}

We now present the main result of this section, that is the characterization of realizable subsets of $\mathbb{Z}^n$.

\begin{theorem}\label{theo:realizable}
A subset $S\subset \mathbb{Z}^n$ is realizable if and only if the following properties hold.
\begin{enumerate}
    \item If $x\in S$, then $x_i\geq 0$ for every $i\in [n]$. Moreover, $x$ has at most one coordinate equal to zero.
    \item For every $i\in [n]$, there exists exactly one element $x\in S$ such that $x_i=0$.
    \item If $x\in S$ and $x_i > 0$, for any $i\in [n]$, then there exists $y\in S$ satisfying $y_i=x_i-1$ and $\max_{j\in [n]} \{\vert y_j-x_j \vert \}\leq 1$.
\end{enumerate}
\end{theorem}

\begin{proof}
Suppose first that the graph $(G,W)$ is a realization of $S$, where $W=\{ \omega_1, \dots, \omega_n\}$. Notice that properties 1 and 2 have been proven in Proposition~\ref{pro:basic}.
To prove property $3$, let $x\in S$ be such that $x_i>0$. Then, $d(u(x),\omega_i)=x_i>0$.
Let $v$ be the neighbor of $u(x)$ in a shortest path from $u(x)$ to $\omega_i$ and let $y=r(v\vert W)$.
Then, $y=r(v\vert W)\in S$, $y_i=d(v,\omega_i)=x_i-1$ and $\vert y_j-x_j\vert =\vert  d(v,\omega_j)-d(u(x),\omega_j)\vert \le 1$, for every $j\in [n]$,
because $u(x)$ and $v$ are neighbors in $G$. Hence,
$\max_{j\in [n]} \{\vert y_j-x_j \vert \}\leq 1$.

Conversely, let $S\subset \mathbb{Z}^n$ be a subset satisfying properties $1$, $2$ and $3$. We define a graph $G=(V,E)$ as follows: $V=S$ and, for $x,y\in S$, $xy\in E$ if and only if  $\max_{j\in [n]} \{\vert y_j-x_j\vert \}=1$.

Then, the set  $W=\{x\in S\colon x_i=0, { \rm  \ for\ some\ } i\in [n]\}$ clearly satisfies $\vert W \vert=n$, by properties $1$ and $2$. We denote by $\omega_i$ the only vertex of $W$ with $i$-th coordinate equal to 0.

Let us see that $(G,W)$ is a realization of $S$. To this end, it is enough to prove that $r(x\vert W)=x$ for every $x\in S$ or, equivalently, $d(x,\omega_i)=x_i$, for every $x\in S$ and $i\in [n]$. Let $x\in S$ and let $i\in [n]$. If $x_i=0$ then, $x=\omega_i$ and $d(x,\omega_i)=d(\omega_i,\omega_i)=0=x _i$.

Now assume that $x_i>0$. Let $\omega_i=v_0,v_1,\dots , v_{k}=x$ be a shortest path between $\omega_i$ and $x$.
Let us see first that $d(x,\omega_i)\geq x_i$. The $i$-th coordinates of two consecutive vertices of $G$ differ in at most one unit, by definition of the edges of $G$. Therefore, $d(x,\omega_i)=k\geq x_i$ because the $i$-th coordinate of $\omega_i$ is zero and the $i$-th coordinate of $x$ is $x_i$.

Now, we prove that $d(x,\omega_i)\leq x_i$. By property 3, there exists $y^1\in S$ such that $y^1_i=x_i-1$ and $\max_{j\in [n]} \{\vert y^{1}_{j}-x_j\vert \}\leq 1$. Therefore, $\max_{j\in [n]} \{\vert y^1_j-x_j \vert \}=\vert y^1_i-x_i\vert =1$ and, consequently, $x,y^1$ are neighbors in $G$. If $y^1_i=x_i-1>0$, by using the same argument as before, there exists $y^2\in S$ such that $y^2_i=y^1_i-1=x_i-2$,  $\max_{j\in [n]} \{\vert y^2_j-y^1_j\vert \}=\vert y^2_i-y^1_i\vert =1$ and $y^2,y^1$ are neighbors in $G$. By repeatedly applying the same argument, we obtain a vertex sequence $x=y^0, y^1, y^2, \dots ,y^{x_i}$ such that $y^j_i=x_i-j$ and $y^j, y^{j+1}$ are neighbors in $G$,
if $0\le j< x_i$. In particular, $y^{x_i}_i=x_i-x_i=0$ and therefore,
    $y^{x_i}=\omega_i$. Finally, since $y^0, y^1, y^2, \dots ,y^{x_i}$ defines a  path of length $x_i$ from $x$ to $\omega_i$, we deduce $d(x,\omega_i)\le x_i$.
\end{proof}

Conditions 1 and 2 in the theorem above are basic properties related with the fact that the coordinates of the metric representation vectors are distances. The meaning of the third condition is also clear, since it allows us to ensure that each vertex at a positive distance, say $d$, from a vertex $w$ of the resolving set, has at least a neighbor at distance $d-1$ from $w$.

\section{Uniqueness}\label{Sec:3}

Once we have characterized the sets of $\mathbb{Z}^n$ that are realizable as the metric representation vector set of a graph, our purpose is now to study the uniqueness of its realizations. In order to put this problem in context, we first show some examples.

\begin{example}\label{example:cycle}
It is straightforward to check that the set $S=\{ (0,3),(3,0), (1,2),$ $ (2,1), (1,4),  (4,1), (2,5), (5,2), (3,4), (4,3)\}$ is realizable by
$(C_{10},W)$ and by $(C_{10}, W')$, if  $V(C_{10})=\{ u_i\colon 1\leq i\leq 10 \}$, $E(C_{10})=\{u_iu_{i+1}\colon 1\leq i< 10\}\cup \{ u_{10} u_1\}$, and $W=\{u_1, u_8\}$, 
$W'=\{u_3, u_6\}$ (see Figure~\ref{fig:cycle_ten}).

\begin{figure}[ht]
     \centering
     \begin{subfigure}{0.35\textwidth}
     \centering
         \includegraphics[width=\textwidth]{cycle_ten_1}
     \end{subfigure}
     \hspace{1.5cm}
     \begin{subfigure}{0.35\textwidth}
     \centering
         \includegraphics[width=\textwidth]{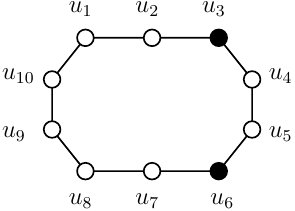}
     \end{subfigure}
        \caption{$(C_{10},\{u_1,u_8\})$ (left) and  $(C_{10},\{u_3,u_6\})$ (right) are realizations of the set $S=\{ (0,3),(3,0), (1,2), (2,1), $ $(1,4), (4,1), $ $(2,5), (5,2), (3,4), (4,3)\}$.}
        \label{fig:cycle_ten}
\end{figure}
\end{example}

\begin{example}\label{example:nonequiv}
In Figure~\ref{fig:two_resolving_sets_a}, $(G, W)$, where $W=\{\omega_1, \omega_2\}$, is a realization of the set $S=\{ (0,3),(3,0),
(1,2),(2,1),  (2,3), (3,2), (1,4), (4,1), (2,5), (5,2),(3,4),  (4,3),  $ $(4,5) ,  (5,4)\}$. Notice that $r(a\vert W)=(2,3)$ and $r(b\vert W)=(3,2)$. In Figure~\ref{fig:two_resolving_sets_b}, $(G, W')$ is a realization of $S$, where $W'=\{z_1, z_2\}$ and, in this case, $r(x\vert W')=(2,3)$ and $r(y\vert W')=(3,2)$.

\begin{figure}[ht]
     \centering
     \begin{subfigure}{.23\textwidth}
     \centering
         \includegraphics[width=\textwidth]{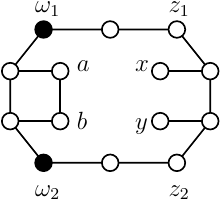}
         \caption{}
         \label{fig:two_resolving_sets_a}
     \end{subfigure}
     \hspace{1.5cm}
     \begin{subfigure}{.23\textwidth}
     \centering
         \includegraphics[width=\textwidth]{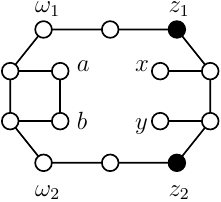}
         \caption{}
         \label{fig:two_resolving_sets_b}
     \end{subfigure}
        \caption{(a) $(G,\{w_1,w_2\})$ and  (b) $(G,\{z_1,z_2\})$ are realizations of the set $S=\{ (0,3),(3,0),
(1,2), (2,1),(2,3), (3,2), (1,4), (4,1), (2,5), (5,2), $ $(3,4), (4,3),  (4,5) ,  (5,4)\}$. The vertices corresponding to the vectors $(2,3)$ and $(3,2)$ are $a$ and $b$ (a), and $x$ and $y$ (b).}
        \label{fig:two_resolving_sets}
\end{figure}
\end{example}

Taking a look to both examples, we notice that the realizations given in Example~\ref{example:cycle} are essentially equal, but this is not the case in Example~\ref{example:nonequiv}, because
the vertices corresponding to the vectors $(2,3)$ and $(3,2)$ are neighbors if we consider the realization $(G,W)$, but they are not when we consider the realization $(G,W')$.

These examples lead us to the following definition.

\begin{definition}\label{def:equivalent}
Let $(G,W)$ and $(G',W')$ be two realizations of a set $S\subset \mathbb{Z}^n$.
We say that both realizations are \emph{equivalent} if the map $f\colon V(G)\to V(G')$ such that $r(u\vert W)=r(f(u)\vert W')$ is a graph isomorphism.
\end{definition}

With this definition, the realizations given in Example~\ref{example:cycle} are equivalent, but not those given in Example~\ref{example:nonequiv}.

Notice that two equivalent realizations of $S\subset \mathbb{Z}^n$ necessarily involve  isomorphic graphs. However, the converse is not true, as shown in  Example~\ref{example:nonequiv} (see Figure~\ref{fig:two_resolving_sets}), where two non equivalent realizations of a set $S$ are given, both sharing the same graph.

Definition~\ref{def:equivalent} provides a suitable framework to study the problem of the uniqueness of the  realization of a set of vectors.

\begin{definition}\label{def:uniquely}
A finite subset $S\subset \mathbb{Z}^n$ is \emph{uniquely realizable} if all of the realizations of $S$ are equivalent.
\end{definition}

\begin{remark}\label{rem:equivrel}
Observe that the definition of equivalent realizations establishes an equivalence relation in the set of all realizations of a fixed set $S\subseteq \mathbb{Z}^n$.
With this terminology, the uniqueness of the realizations consists of characterizing all sets $S$ such that there is only  one equivalence class.
\end{remark}

In order to address the problem of the uniqueness of the realization, notice that for every realizable set of vectors $S$, it is possible to build a graph $G$ and a resolving set $W$ realizing $S$ with the procedure described in the proof of Theorem~\ref{theo:realizable}. For this reason, we propose the following definition.

\begin{definition}\label{def:canonical}
Let $S\subset \mathbb{Z}^n$ be a realizable set. The \emph{canonical realization} of $S$ is $(\widehat{G},\widehat{W})$,
where
$$V(\widehat{G})=S, \ \ E(\widehat{G})=\{ xy\colon x,y\in S \textrm{ and } \max_{i\in [n]} \{\vert y_i-x_i\vert \}=1\},$$
and
$$\widehat{W}=\{x\in S\colon x_i=0, {\rm \ for\ some\ } i\in [n]\}.$$
\end{definition}

\begin{remark}
Since the canonical realization exists for every realizable set, a set $S$ of vectors is uniquely realizable if and only if all the realizations of $S$ are equivalent to the canonical one.
\end{remark}

Thus, we now focus on studying the structure of the canonical realization, that can be viewed as a subgraph of the strong product of paths.

The strong product of graphs was introduced in~\cite{Sabidussi1959} and extensive information on this concept can be found in~\cite{Imrich2000}.
The strong product of two graphs $G$ and $H$, denoted by $G\boxtimes H$, is the graph such that $V(G\boxtimes H)=V(G)\times V(H)$ and two vertices $(g,h), (g',h')\in V(G)\times V(H)$ are adjacent in $G\boxtimes H$ if and only if one of the following conditions holds:
\begin{itemize}
    \item $g=g'$ and $h,h'$ are adjacent in $H$,
    \item $h=h'$ and $g,g'$ are adjacent in $G$,
    \item $g,g'$ are adjacent in $G$ and $h,h'$ are adjacent in $H$.
\end{itemize}
The distance between two vertices in $G\boxtimes H$ is given by
$$d((g,h), (g',h'))=\max \{d_G(g,g'), d_H(h,h')\}.$$
The natural extension of this operation to a finite family of graphs $\{G_i\colon i\in[n]\}$ provides the graph $G_1\boxtimes \dots \boxtimes G_n$ with vertex set $V(G_1)\times \dots \times V(G_n)$ and such that two vertices $(g_1,\dots ,g_n)$ and $(g_1',\dots,g_n')$ from $ V(G_1)\times \dots \times V(G_n)$ are adjacent if and only if $\max_{i\in [n]}\{ d_{G_i}( g_i,g'_i) \}=1$.

Let $P_r$ denote the path graph of order $r$ with vertex set $V(P_r)=\{0,1,\dots , r-1\}$, where
$i-1$ and $i$ are adjacent for every $i\in [r-1]$.
Therefore,
$d_{P_r}(i,j)=\vert i-j\vert$ for any two vertices $i, j$ in $P_r$.
Following this notation, the strong product of $n$ copies of $P_r$ is the graph $P_r \boxtimes \stackrel{n)}{\dots} \boxtimes P_r$, with vertex set $\{0,1,\dots, r-1\}^n$ and two vertices $x, y$ are  adjacent if and only if $\max_{i\in [n]}\{ \vert x_i-y_i\vert \}=1$.

In addition, the  graph $\widehat{G}$ described in Definition~\ref{def:canonical} is clearly the subgraph induced by $S$ in the strong product of $n$ paths of order at least $1+\max \{x_i\colon x\in S, i\in [n]\}$.
We next prove that every graph $G$ of a realization of a subset $S$ of $\mathbb{Z}^n$ is isomorphic to a subgraph of the strong product of $n$ paths, not necessarily induced.

\begin{proposition}\label{pro:isomorphic}
Let $(G,W)$ be a realization of a set $S\subset \mathbb{Z}^n$.
The map $\phi\colon V(G)\to S$ such that $\phi(v)=r(v\vert W)$ defines a graph isomorphism between $G$ and a subgraph of the strong product of $n$ paths of order $r$, for any $r\ge 1+\max \{x_i\colon x\in S, i\in [n]\}$.
\end{proposition}

\begin{proof}
The choice of $r$ and Condition $1$ of Theorem~\ref{theo:realizable} ensure that the map is well defined. Moreover, $\phi$ is injective, because $W$ is a resolving set, and $\phi (V)=S$, because $(G,W)$ is a realization of $S$. Hence, $\phi$ is a bijection.

Now, we prove that $\phi(v)\phi(v')\in E(P_r \boxtimes  \stackrel{n)}{\dots}  \boxtimes P_r)$, if $vv'\in E(G)$. Indeed, recall that $v=u(x)$ and $v'=u(y)$ for some $x,y\in S$, such that $\phi(v)=r(v\vert W)=x$ and $\phi(v')=r(v'\vert W)=y$.  Hence,  $\vert y_i-x_i\vert=\vert d(u(y),\omega_i)-d(u(x),\omega_i)\vert=\vert d(v',\omega_i)-d(v,\omega_i)\vert \leq 1$, for every $i\in [n]$, because $v$ and $v'$ are neighbors in $G$. Therefore, $\max_{i\in [n]}\{ \vert y_i-x_i\vert \}\leq 1$.
Moreover, $x=r(v\vert W)\neq r(v'\vert W)=y$, because $v\neq v'$. Thus, there exists $j\in [n]$ such that $x_j\neq y_j$, implying that $\max_{i\in [n]}\{ \vert y_i-x_i\vert \}=1$. This means that $xy\in E(P_r \boxtimes  \stackrel{n)}{\dots}  \boxtimes P_r)$, where $x=\phi(v)$ and $y=\phi(v')$. Hence, $G$ is isomorphic to the graph $G(S,\{\phi (v)\phi (v' ) : vv'\in E(G)\})$, that is a subgraph of the strong product of $n$ paths of order at least $1+\max \{x_i\colon x\in S, i\in [n]\}$.
\end{proof}

A first consequence of this construction is that we can provide a generalization to every dimension of the property suggested  in~\cite{Khuller1996} and shown in~\cite{Claverol2021}, stating that every graph with metric dimension at most two can be represented as a  subgraph of the strong product of two paths. Notice that this subgraph is non necessarily induced.

\begin{corollary}\label{cor:embedding}
Let $G$ be a graph with $\dim (G)\leq n$ and let $r$ be an integer such that
$r\geq 1+diam(G)$.
Then, $G$ is isomorphic to a subgraph of $P_r \boxtimes \stackrel{n)}{\dots}  \boxtimes P_r$, not necessarily induced.
\end{corollary}

It was proved in~\cite{Graham1985} that any graph can be canonically isometrically embedded into the Cartesian product of graphs. Following this idea, the isometric dimension of a graph $G$ is the least number of graphs needed to isometrically embed $G$ into a Cartesian product. Moreover, the strong isometric dimension of a graph~\cite{Fitzpatrick2000,Jerebic2006} is the least number $k$ such that $G$ embeds isometrically into the strong product of $k$ paths. A general reference about isometric embeddings of graphs into product graphs can be found in~\cite{Imrich2000}.

Corollary~\ref{cor:embedding} has a similar flavor but regarding the metric dimension. However, notice that the isomorphism we are using is not an isometric embedding, because it just preserves the distances from the vertices in the resolving set to any other vertex of the graph, but not necessarily between every pair of vertices.
\vspace{3mm}

The map $\phi$ defined in Proposition~\ref{pro:isomorphic} provides a useful tool to decide whether two realizations  of a set of $\mathbb{Z}^n$ are equivalent, since all realizations of $S$ can be seen as subgraphs of a common graph.
Moreover, since the subgraphs of the strong product of paths associated to the realizations of a fixed $S$ share the same vertex set, the equivalency between two realizations can be characterized by means of the edge sets.
We will use the subgraph of the strong product of graphs considered in Proposition~\ref{pro:isomorphic} to explore in depth such relationship.

\begin{definition}\label{def:inmersion}
Let $(G,W)$ be a realization of $S\subset \mathbb{Z}^n$. We denote by $G^*_W$ the subgraph of the strong product of $n$ paths of order $r$,  with $r\geq 1+diam(G)$, whose vertex set is $V(G^*_W)=S$ and whose edge set is $E(G^*_W)=\{ xy : x=r(v\vert W), y=r(v'\vert W),vv'\in E(G)\}$.
\end{definition}

We illustrate this definition with some examples.

\begin{example}\label{example:immersion_cycle}
In Figure~\ref{fig:inmmersion_cycle_a} we show the cycle with ten vertices $C_{10}$ and its resolving set $W$, the set of black vertices. Note that $(C_{10},W)$ is a realization of $S=\{ (0,3),(3,0), (1,2),(2,1), (1,4), (4,1), (2,5), (5,2), (3,4), (4,3)\}$. The subgraph $G^*_W$ of $P_8\boxtimes P_8$ defined above, is shown in Figure~\ref{fig:inmersion_cycle_b}.

\begin{figure}[h]
     \centering
     \begin{subfigure}[b]{.3\textwidth}
     \centering
         \includegraphics[width=\textwidth]{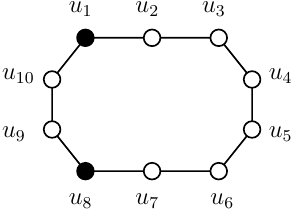}
         \caption{}
         \label{fig:inmmersion_cycle_a}
     \end{subfigure}
     \hspace{2cm}
     \begin{subfigure}[b]{.3\textwidth}
     \centering
         \includegraphics[width=\textwidth]{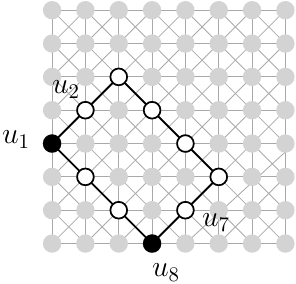}
         \caption{}
         \label{fig:inmersion_cycle_b}
     \end{subfigure}

        \caption{(a) A realization $(G,W)$  of the set $S=\{ (0,3),(3,0),
(1,2), (2,1), $ $(1,4), (4,1), (2,5), (5,2),(3,4), (4,3)\}$ and (b) the graph $G_W^*$.}
        \label{fig:inmersion_cycle}
\end{figure}
\end{example}

\begin{example}\label{example:immersion}
In Figure~\ref{fig:inmmersion_a} we show a graph $G$ and its resolving set $W$, the set of black vertices. Note that $(G,W)$ is a realization of $S=\{ (0,3),(3,0), $ $ (1,2), (2,1), (2,3),  (3,2), (1,4), (4,1), (2,5), (5,2), (3,4), (4,3), (4,5), $ $ (5,4)\}$. The subgraph $G^*_W$ of $P_6\boxtimes P_6$ defined above, is shown in Figure~\ref{fig:inmersion_b}.
\begin{figure}[h]
     \centering
     \begin{subfigure}[b]{.25\textwidth}
     \centering
         \includegraphics[width=\textwidth]{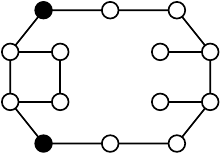}
         \caption{}
         \label{fig:inmmersion_a}
     \end{subfigure}
     \hspace{2.3cm}
     \begin{subfigure}[b]{.25\textwidth}
     \centering
         \includegraphics[width=\textwidth]{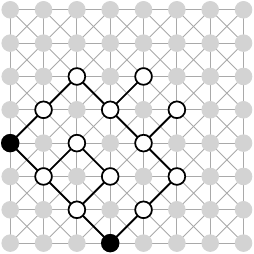}
         \caption{}
         \label{fig:inmersion_b}
     \end{subfigure}

        \caption{(a) A realization $(G,W)$  of the set $S=\{ (0,3),(3,0),
(1,2), $ $(2,1), $ $(2,3),  $ $(3,2),  $ $(1,4), (4,1), (2,5), (5,2), (3,4), (4,3),  (4,5) ,  (5,4)\}$ and (b) the graph $G_W^*$.}
        \label{fig:inmersion}
\end{figure}
\end{example}

Observe that in the graph $G^*_W$, the distance from the vertex $w_i\in W$ to any other vertex $x\in S$ is precisely the $i$-th coordinate of $x$, $x_i$. Moreover, the $i$-th coordinates of the vertices of any shortest path from $w_i$ to $x$ are $0,1,2,\dots , x_i$, respectively.

Notice also that, as we have already pointed out, the subgraph $G^*_W$ is not necessarily an induced subgraph of the strong product of paths. Moreover, by definition, $G^*_W$ preserves the distances from a vertex in $W$ to any other vertex, but not necessarily between two any vertices of $G^*_W$, even if it is an induced subgraph. For instance, in Example~\ref{example:immersion_cycle} (see Figure~\ref{fig:inmersion_cycle}), $G^*_W$ is an induced subgraph of the strong product of paths, however, the vertices $u_2$ and $u_7$ are at distance $5$ in $G$, but the corresponding vertices in the strong product of paths, $(1,4)$ and $(4,1)$ are at distance 3.

The following result is a direct consequence of Proposition~\ref{pro:isomorphic}.

\begin{corollary}
\label{cor:inmersion}
Let $(G,W)$ be a realization of $S\subset \mathbb{Z}^n$. Then, $(G^*_W,\phi(W))$ is a realization of $S$  equivalent to $(G,W)$.
\end{corollary}

\begin{corollary}\label{cor:equality}
Let $(G,W)$ and $(G',W')$ be two realizations of a set $S\subset \mathbb{Z}^n$. Then, $(G,W)$ and $(G',W')$ are equivalent if and only if $E(G^*_W)=E(G'^*_{W'})$.
\end{corollary}

\begin{proof}
Assume that $G$ and $G'$ are equivalent realizations of $S$ and denote by  $f\colon V(G)\to V(G')$ the map such that $r(u\vert W)=r(f(u)\vert W')$, that is a graph isomorphism, by hypothesis. Therefore, $x=r(u(x)\vert W)=r(f(u(x))\vert W')$. We know that $V(G^*_W)=V(G'^*_{W'})=S$, and we have to prove that $E(G^*_W)=E(G'^*_{W'})$. Let $x,y\in S$ then, $xy\in E(G^*_W)$ if and only if $u(x)u(y)\in E(G)$ if and only if $f(u(x))f(u(y))\in E(G')$ if and only if $xy\in E(G'^*_{W'})$.

Conversely, if $E(G^*_W)=E(G'^*_{W'})$ them, both subgraphs of $P_r \boxtimes  \stackrel{n)}{\dots}  \boxtimes P_r$ are equal. Finally, both $G$ and $G'$ are equivalent to the same subgraph of $P_r \boxtimes  \stackrel{n)}{\dots}  \boxtimes P_r$, so we obtain that $G$ and $G'$ are equivalent realizations of $S$.
\end{proof}

\begin{corollary}
    A set $S$ is uniquely realizable if and only if $E(G_W^*)=E(\widehat{G})$ for every realization $(G,W)$ of $S$.
\end{corollary}

In the following example we show the subgraphs of the strong product equivalent to three realizations of the same set $S$, including the canonical one.

\begin{example}
We show three subgraphs of $P_6\boxtimes P_6$ realizing the set $S=\{ (0,3),(3,0),  $ $ (1,2), (2,1), (2,3), (3,2), (1,4),(4,1), (2,5), (5,2),(3,4), (4,3), (4,5), (5,4)\}$. In each case, subgraph edges are the black ones, the vertex set are the black and the white vertices and the black vertices are the resolving set. The vertex set is $S$ in all cases.

In Figure~\ref{fig:strong_product_a} and Figure~\ref{fig:strong_product_b} we show the realizations given in Figure~\ref{fig:two_resolving_sets_a} and Figure~\ref{fig:two_resolving_sets_b} respectively,  and as subgraphs of $P_6\boxtimes P_6$. We can easily check that these subgraphs of $P_6\boxtimes P_6$ have different edge sets, so they are non-equivalent realizations of the set $S$. Moreover, none of them is an induced subgraph, that means that they are not the canonical realization of $S$. Finally, the canonical realization of $S$ is the induced subgraph of $P_6\boxtimes P_6$ shown in Figure~\ref{fig:strong_product_c}.

\begin{figure}[h]
     \centering
     \begin{subfigure}{.2\textwidth}
     \centering
         \includegraphics[width=\textwidth]{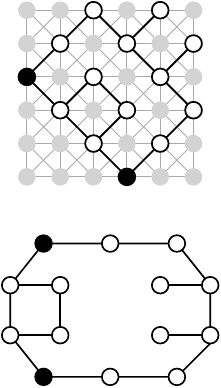}
         \caption{ }
         \label{fig:strong_product_a}
     \end{subfigure}
     \hspace{1cm}
     \begin{subfigure}{.2\textwidth}
     \centering
         \includegraphics[width=\textwidth]{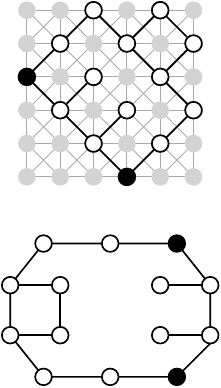}
         \caption{\ }
         \label{fig:strong_product_b}
     \end{subfigure}
    \hspace{1cm}
     \begin{subfigure}{.2\textwidth}
     \centering
         \includegraphics[width=\textwidth]{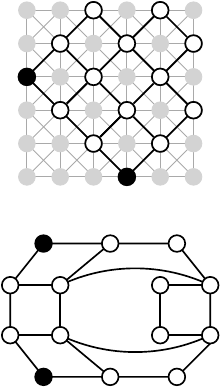}
         \caption{\ }
         \label{fig:strong_product_c}
     \end{subfigure}
        \caption{Three non equivalent realizations of a set $S$ as subgraphs of $P_6\boxtimes P_6$.}
        \label{fig:strong_product}
\end{figure}
\end{example}

It is clear that not every subgraph of $P_r \boxtimes \stackrel{n)}{\dots}  \boxtimes P_r$ with vertex set $S$ is a realization of such set. However, we show in the following result that adding any available edge of the strong product to a subgraph realizing $S$ provides a non equivalent realization of it.

\begin{proposition}\label{prop:adding_edges}
Let $G$ be a subgraph of a strong product of $n$ paths such that $(G,W)$ is a realization of a subset $S\subset \mathbb{Z}^n$ not equivalent to the canonical one. Consider a graph $G'$ obtained by adding to $G$ an edge of the strong product not belonging to $G$. Then, $(G',W)$ is a  realization of $S$ not equivalent to $(G,W)$.
\end{proposition}
\begin{proof}
Let $x,y\in V(G)$ be such that $x$ and $y$ are neighbors in the strong product of paths, but not in $G$. Consider the graph $G'=(V(G),E(G)\cup \{ e\})$, where $e=xy$. Let $z\in V(G')$. Notice that $d_{G}(z,\omega_i)=z_i=d_{P_r \boxtimes  \stackrel{n)}{\dots}  \boxtimes P_r} (z,\omega_i)$ and therefore, the addition of an edge of the strong product $P_r \boxtimes \stackrel{n)}{\dots}  \boxtimes P_r$ cannot decrease such distance. Thus, $d_{G'}(z,\omega_i)=d_{G}(z,\omega_i)=z_i$. Hence, $W$ is a resolving set of $G'$ such that $\{r(z\vert W)\colon z\in V(G')\}=S$.

Finally, it is clear that $G$ and $G'$ have different edge sets, so by Corollary~\ref{cor:equality}, $(G,W)$ and $(G',W)$ are not equivalent realizations.
\end{proof}
\par\medskip

\section{Uniqueness of the realization in the case $n=2$}\label{subsec:case_2}

We now focus on the case $n=2$, in order to characterize the vector sets which are uniquely realizable. We know that every realizable set $S\subset \mathbb{Z}^2$ admits the canonical realization that is an induced subgraph of $P_r\boxtimes P_r$. Therefore, this realization will be the only one if it is not possible to obtain another one by removing edges from it.

Contrary to what happens with the addition of edges (see Proposition~\ref{prop:adding_edges}), it is clear that removing edges from a subgraph of $P_r \boxtimes  \stackrel{n)}{\dots}  \boxtimes P_r$ and keeping the distances between its vertices is not always possible. We focus in these ideas to characterize the uniquely realizable vector sets.

We begin by proving some technical lemmas.

\begin{lemma}\label{lem:first_coordinate}
Let $(\widehat{G},\widehat{W})$ be the canonical realization of $S\subset \mathbb{Z}^2$, with $\widehat{W}=\{\omega_1, \omega_2\}$. Suppose that $x\in S$ and there is  a shortest path in $\widehat{G}$ between $\omega_1$ and $x$ going through an edge $uv$. If at least one of the following conditions holds for some integers $\alpha$ and $\beta$:
\begin{enumerate}[i)]
    \item $u=(\alpha,\beta+1), v=(\alpha+1,\beta), t=(\alpha,\beta)\in S$;
    \item $u=(\alpha,\beta+1), v=(\alpha+1,\beta+1), t=(\alpha,\beta)\in S$;
    \item $u=(\alpha,\beta+2), v=(\alpha+1,\beta+1), t=(\alpha,\beta)\in S$;
    \item $u=(\alpha,\beta), v=(\alpha+1,\beta+1), t=(\alpha,\beta+2)\in S$;
\end{enumerate}
then, $\omega_1\neq u$ and there exists a shortest path between $\omega_1$ and $x$ avoiding the edge $uv$. Moreover, in such a case, the graph $\widehat{G}-uv$ is connected.
\end{lemma}

\begin{proof}
First, in all cases $u\neq t$ and the first coordinate of $u$ equals the first coordinate of $t$, so that this coordinate cannot be zero, by property 2 of Theorem~\ref{theo:realizable}. This means that $u\neq \omega_1$. Moreover, the equality of the first coordinates implies $d(u,\omega_1)=d(t,\omega_1)$ in all cases. Notice also that $t$ and $v$ are neighbors in $\widehat{G}$, by definition of $\widehat{G}$. Therefore, we can construct a shortest path between $\omega_1$ and $x$ by joining a shortest path from $w_1$ to $t$, the edge $tv$ and a shortest path from $v$ to $x$, that clearly avoids the edge $uv$ (see Figure~\ref{fig:lemafirstcoordinates}).
\begin{figure}[h]
    \centering
    \includegraphics[width=0.65\textwidth]{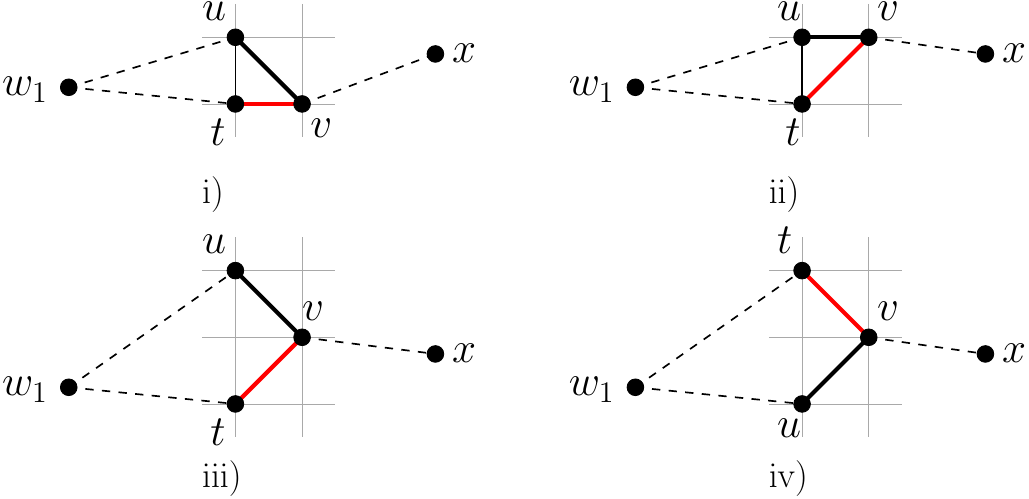}
    \caption{In all cases, if there is a shortest path from $w_1$ to a vertex $x$ that goes through $uv$, then there is a shortest path from $w_1$ to $x$  that goes through $tv$ avoiding $uv$.}
    \label{fig:lemafirstcoordinates}
\end{figure}

Now, let $x,y\in S$. We have seen that there are paths in $\widehat {G}$ from $w_1$ to $x$ and to $y$  avoiding the edge $uv$. Hence, it is possible to build a path in $\widehat {G}$ between $x$ and $y$ avoiding the edge $uv$. Hence, $\widehat{G}-uv$ is connected.
\end{proof}

By interchanging the role of both coordinates in the preceding lemma, we obtain the following result.

\begin{lemma}\label{lem:second_coordinate}
Let $(\widehat{G},\widehat{W})$ be the canonical realization of $S\subset \mathbb{Z}^2$, with $\widehat{W}=\{\omega_1, \omega_2\}$. Suppose that $x\in S$ and there is  a shortest path in $\widehat{G}$ between $\omega_2$ and $x$ going through an edge $uv$. If one of the following conditions hold:
\begin{enumerate}[i)]
    \item $u=(\alpha+1,\beta), v=(\alpha,\beta+1), t=(\alpha,\beta)\in S$;
    \item $u=(\alpha+1,\beta), v=(\alpha+1,\beta+1), t=(\alpha,\beta)\in S$;
    \item $u=(\alpha+2,\beta), v=(\alpha+1,\beta+1), t=(\alpha,\beta)\in S$;
    \item $u=(\alpha,\beta), v=(\alpha+1,\beta+1), t=(\alpha+2,\beta)\in S$;
\end{enumerate}
then, $\omega_2\neq u$ and there exists a shortest path between $\omega_2$ and $x$ avoiding the edge $uv$. Moreover, in such a case, the graph $\widehat{G}-uv$ is connected.
\end{lemma}

\begin{figure}[h]
    \centering
    \includegraphics[width=0.25\textwidth]{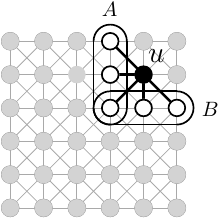}
    \caption{At least one vertex in $A$ and at least one vertex in $B$ belong to $S$ and are adjacent to $u$.}
    \label{fig:neighbors}
\end{figure}

\begin{lemma}\label{lem:neighbors_AB}
Let $(G,W)$ be a realization of $S\subset \mathbb{Z}^2$, with $G$ a subgraph of the strong product $P_r \boxtimes   P_r$.
Let $u=(u_1,u_2)\in S$. If $u_1>0$, then at least one vertex in the set $A=\{(u_1-1, u_2-1), (u_1-1, u_2), (u_1-1,u_2+1)\}$ belongs to $S$ and is adjacent to $u$. If $u_2>0$, then at least one vertex in $B=\{(u_1-1, u_2-1), (u_1, u_2-1), (u_1+1,u_2-1)\}$ belongs to $S$ and is adjacent to $u$. (See Figure~\ref{fig:neighbors}).
\end{lemma}
 \begin{proof}
Suppose that $W=\{w_1,w_2\}$ and $u=(u_1,u_2)=(d(w_1,u), d(w_2,u))$. If $u_1>0$, then there is a neighbor of $u$ in $S$ at distance $u_1-1$ from $w_1$. Hence, at least one of the vertices of $A$ belongs to $S$ and is adjacent to $u$. Similarly, if $u_2>0$, we derive that one of the vertices of $B$ belongs to $S$ and is adjacent to $u$.
\end{proof}
We now prove the main result of this section that characterizes the subsets of $\mathbb{Z}^2$ that are uniquely realizable. They are characterized in terms of conditions of the coordinates of their elements.

\begin{theorem}\label{theo:uniqueness}
Let $S\subset \mathbb{Z}^2$ be a realizable set. Then, there exist two non equivalent realizations of $S$ if an only if at least one of the following conditions holds, for some integers $\alpha,\beta$:
\begin{enumerate}
    \item $\{(\alpha,\beta),(\alpha,\beta+1),(\alpha+1,\beta)\}\subset S$ (Figure~\ref{fig:case1}),
    \item $\{(\alpha,\beta),(\alpha,\beta+1),(\alpha+1,\beta+1)\}\subset S$ (Figure~\ref{fig:case2}),
    \item $\{(\alpha,\beta),(\alpha+1,\beta),(\alpha+1,\beta+1)\}\subset S$ (Figure~\ref{fig:case3}),
    \item $\{(\alpha-1,\beta),(\alpha,\beta-1),(\alpha,\beta+1),(\alpha+1,\beta)\}\subset S$ (Figure~\ref{fig:case4})
    \item $\{(\alpha,\beta),(\alpha,\beta+2),(\alpha+1,\beta+1),(\alpha+2,\beta)\}\subset S$ (Figure~\ref{fig:case5}).
\end{enumerate}

\begin{figure}[h]
     \centering
     \begin{subfigure}{.20\textwidth}
     \centering
         \includegraphics[width=\textwidth]{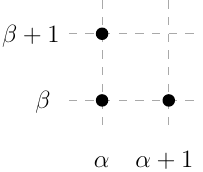}
         \caption{ }
         \label{fig:case1}
     \end{subfigure}
     \hspace{0.8cm}
     \begin{subfigure}{.20\textwidth}
     \centering
         \includegraphics[width=\textwidth]{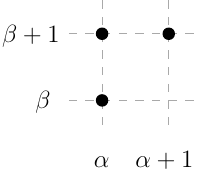}
         \caption{\ }
         \label{fig:case2}
     \end{subfigure}
    \hspace{0.8cm}
     \begin{subfigure}{.20\textwidth}
     \centering
         \includegraphics[width=\textwidth]{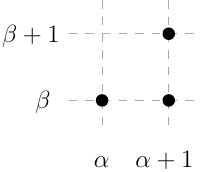}
         \caption{\ }
         \label{fig:case3}
     \end{subfigure}
      \hspace{0.8cm}
     \begin{subfigure}{.26\textwidth}
     \centering
         \includegraphics[width=\textwidth]{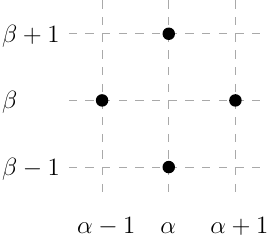}
         \caption{\ }
         \label{fig:case4}
     \end{subfigure}
      \hspace{0.8cm}
     \begin{subfigure}{.26\textwidth}
     \centering
         \includegraphics[width=\textwidth]{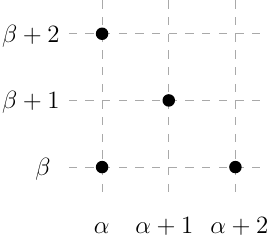}
         \caption{\ }
         \label{fig:case5}
     \end{subfigure}
        \caption{The five configurations of vertices in $S\subset \mathbb{Z}^2$ of  Theorem~\ref{theo:uniqueness}.}
        \label{fig:uniqueness}
\end{figure}

\end{theorem}

\begin{proof}

We first proof the necessity. Let $(G,W)$ and $(H,R)$ be two non-equivalent realizations of $S$. By Corollaries~\ref{cor:inmersion} and \ref{cor:equality}, we may assume that $G$ and $H$ are subgraphs of a strong product $P_r  \boxtimes P_r$ satisfying $V(G)=V(H)=S$ and $E(G)\neq E(H)$. Without lost of generality, we may also assume that there exist two vertices $x$ and $y$ in $S$ that are neighbors in $H$, but not in $G$. In particular, these vertices are neighbors in $P_r \boxtimes P_r$, so that their coordinates differ in at most one unit.
We distinguish the following cases.

\begin{enumerate}[i)]
    \item The vertices $x$ and $y$ have the same first coordinate.
    Assume that $x=(a,b+1)$ and $y=(a,b)$. In particular, $b+1>0$.
    Since $x$ and $y$ are not neighbors in $G$, by applying Lemma~\ref{lem:neighbors_AB} to $x$, we derive that either $z=(a+1,b)\in S$, and then Condition 1 holds for $\{x,y,z\}$,
    or $t=(a-1,b)\in S$, and then Condition 3 holds for $\{x,y,t\}$.

    \item The vertices $x$ and $y$ have the same second coordinate.
    We proceed as in the preceding case but interchanging the role of the first and second coordinates, and we obtain that either Condition 1 or Condition 2 holds.

    \item
    The first and second coordinates of one of the vertices  are less and greater than the first and second coordinates of the other vertex, respectively.
    We may assume that $x=(a,b+1)$ and $y=(a+1,b)$. In particular, $a+1>0$ and $b+1>0$.
    If $z=(a,b)\in S$, then $\{x,y,z\}$ satisfies Condition 1. Otherwise, since $x$ and $y$ are not neighbors in $G$, by applying Lemma~\ref{lem:neighbors_AB} to $x$ for the set $B$, we obtain that $s=(a-1,b)\in S$,
    and by applying Lemma~\ref{lem:neighbors_AB} to $y$ for the set $A$, we obtain that $t=(a,b-1)\in S$. Thus, $\{x,y,s,t\}$ satisfies Condition 4.

    \item
    One vertex has both coordinates greater than those of the other vertex. Then, we may assume $x=(a+1,b+1)$, $y=(a,b)$. In particular, $a+1>0$ and $b+1>0$.

    If $z=(a, b+1)\in S$, then $\{x,y,z\}$ satisfies Condition 2. If $s=(a+1, b)\in S$ then $\{x,y,s\}$ satisfies Condition 3. Now suppose that $z,s\notin S$. Then, since $x$ and $y$ are not neighbors in $G$, by applying Lemma~\ref{lem:neighbors_AB} to $x$ for both sets $A$ and $B$, we obtain that $t=(a, b+2)\in S$ and $p=(a+2, b)\in S$. Therefore $\{x,y,t,p\}$ satisfies Condition 5.
\end{enumerate}

We now prove the sufficiency. Let $(\widehat{G},\widehat{W})$ be the canonical realization of $S$, with $\widehat{W}=\{\omega_1, \omega_2\}$, that exists because $S$ is realizable. Assume that there exist integers $\alpha, \beta$ such that one of the conditions from 1 to 5 hold. To conclude the proof, we will see that there is an edge $e$ of $\widehat{G}$ such that $(\widehat{G}-e, \widehat{W})$ is also a realization of $S$, but not equivalent to the canonical one, because the graphs $\widehat{G}$ and $\widehat{G}-e$ have different size and, consequently, are not isomorphic.
For this purpose, in each one of the five cases we choose an edge $e$ of $\widehat{G}$
such that the deletion of $e$ does not change the distances from the vertices in $W$ to any other vertex.
\begin{enumerate}
    \item If $\{(\alpha,\beta),(\alpha,\beta+1),(\alpha+1,\beta)\}\subset S$, let
    $u=(\alpha,\beta+1)$ and $v=(\alpha+1,\beta)$.

    Since $t=(\alpha,\beta)\in S$, if there is a shortest path in $\widehat{G}$ between $\omega_1$ and a vertex $x$ going through the edge $uv$, then  there exists a shortest path between $\omega_1$ and $x$ avoiding the edge $uv$, because Condition i) of  Lemma~\ref{lem:first_coordinate} is satisfied for $u$, $v$ and $t$.

     Similarly, if there is a shortest path  in $\widehat{G}$ between $\omega_2$ and a vertex $y$ going through the edge $uv$, then there exists a shortest path between $\omega_2$ and $y$ avoiding the edge $uv$, because Condition i) of Lemma~\ref{lem:second_coordinate}$i)$ is satisfied for $u$, $v$ and $t$.

    \item If $\{(\alpha,\beta),(\alpha,\beta+1),(\alpha+1,\beta+1)\}\subset S$,  let $u=(\alpha,\beta+1)$ and $v=(\alpha+1,\beta+1)$.

    Since $t=(\alpha,\beta)\in S$, if there is a shortest path in $\widehat{G}$ from $w_1$ to $x$ going trough $uv$, we derive that there is a shortest path from $\omega_1$ to  $x$ avoiding the edge $uv$, because Condition ii) of Lemma~\ref{lem:first_coordinate} is satisfied for $u$, $v$ and $t$.

    Notice that, in this case, there are no shortest paths between $\omega_2$ and any other vertex $y$ going through the edge $uv$, because $u$ and $v$ have the same second coordinate.

    \item If $\{(\alpha,\beta),(\alpha+1,\beta),(\alpha+1,\beta+1)\}\subset S$,
    let $u=(\alpha+1,\beta)$ and $v=(\alpha+1,\beta+1)$.

    Since $t=(\alpha,\beta)\in S$, if there is a shortest path in $\widehat{G}$  from $w_2$ to $y$ going through the edge $uv$, then we derive that there is also a shortest path from  $\omega_2$ to $y$ avoiding the edge $uv$, because Condition ii) of Lemma~\ref{lem:second_coordinate} is satisfied for $u$, $v$ and $t$.

     Notice that there are no shortest paths in $\widehat{G}$ between $\omega_1$ and any other vertex $x$ going through the edge $uv$, because the first coordinates of $u$ and $v$ are equal.

    \item
    If $\{(\alpha-1,\beta),(\alpha,\beta-1),(\alpha,\beta+1),(\alpha+1,\beta)\}\subset S$, let
    $u=(\alpha,\beta+1)$ and $v=(\alpha+1,\beta)$.

    Since $(\alpha,\beta-1)\in S$, if there is a shortest path in $\widehat{G}$ from $w_1$ to $x$  going through the edge $uv$, then there is also a shortest path from $w_1$ to $x$ avoiding the edge $uv$, because Condition iii) of Lemma~\ref{lem:first_coordinate} is satisfied.

    Since $(\alpha-1,\beta)\in S$, if there is a shortest path in $\widehat{G}$ from $w_2$ to $y$  going through the edge $uv$, then there is also a shortest path from $w_2$ to $y$ avoiding the edge $uv$,
    because Condition iii) of Lemma~\ref{lem:second_coordinate} is satisfied.

    \item If $\{(\alpha,\beta),(\alpha,\beta+2),(\alpha+1,\beta+1),(\alpha+2,\beta)\}\subset S$, let $u=(\alpha,\beta)$ and $v=\alpha+1,\beta+1)$.

    Since $(\alpha,\beta+2)\in S$, if there is a shortest path in $\widehat{G}$ from $w_1$ to $x$  going through the edge $uv$, then there is also a shortest path from $w_1$ to $x$ avoiding the edge $uv$, because Condition iv) of Lemma~\ref{lem:first_coordinate} is satisfied.

    Since $(\alpha+2,\beta)\in S$, if there is a shortest path in $\widehat{G}$ from $w_2$ to $y$  going through the edge $uv$, then there is also a shortest path from $w_2$ to $y$ avoiding the edge $uv$,
    because Condition iv) of Lemma~\ref{lem:second_coordinate} is satisfied.

\end{enumerate}

In all cases, the choice of the edge $e=uv$ ensures that the graph $H=\widehat{G}-e$ is connected, as is pointed out in Lemmas~\ref{lem:first_coordinate} and ~\ref{lem:second_coordinate}. In addition, $d_{H}(z,\omega_i)=d_{\widehat{G}}(z,\omega_i)$, for every $z\in S$ and $i\in \{1,2\}$ and thus, $\widehat{W}$ is a resolving set of $H$ satisfying $\{r(z\vert \widehat{W})\colon z\in V(H)\}=S$. Moreover, $H$ and $\widehat{G}$ are non-isomorphic graphs because they  have different number of edges. Therefore, $(\widehat{G},\widehat{W})$ and $(H,\widehat{W})$ are non-equivalent realizations of $S$.
\end{proof}

\bmhead{Acknowledgements} This work was partially supported by the grants of the Spanish Ministry of Science and Innovation PID2021-123278OB-I00 and PID2019-104129GB-I00, both
funded by MCIN/AEI/ 10.13039/501100011033 and by ``ERDF A way of making Europe'', and Gen. Cat. DGR 2021-SGR-00266.

\bibliographystyle{spmpsci}
\bibliography{bibfile}

\par\bigskip

\bmhead{Declarations}

\bmhead{Conflict of interest} The authors declare no conflict of interest.

\bmhead{Author Contributions} Both authors contributed equally to this work.

\end{document}